\documentclass[11pt,reqno]{amsart}
\usepackage{graphicx,amsmath,amsthm,verbatim,amssymb}
\usepackage{lineno,titletoc}
\usepackage{mathrsfs}
\usepackage{adjustbox}
\usepackage{ytableau}
\usepackage{bbm}
\usepackage{hyperref}
\hypersetup{colorlinks}
\usepackage[numbers]{natbib}
\usepackage{array}

\usepackage[font=scriptsize]{caption}

\usepackage[T1]{fontenc} 
\usepackage{fourier} 
\usepackage[english]{babel} 
\usepackage{amsmath,amsfonts,amsthm}
\usepackage{appendix}

\usepackage[all]{xy}
\usepackage{enumerate}
\usepackage{tikz}
\usepackage{mathrsfs}
\usepackage[english]{babel}
\usepackage{multirow}
\usepackage{float}
\usepackage{enumitem}
\usepackage{graphicx,accents}

\newcolumntype{M}[1]{>{\centering\arraybackslash}m{#1}}
\newcolumntype{N}{@{}m{0pt}@{}}

\pdfoutput=1

\def\Z{\mathbb{Z}}

\def\L{\mathcal{L}}

\def\la{\leftarrow}
\def\ra{\rightarrow}
\def\r{\mathtt{r}}

\def\eps{\varepsilon}

\newtheorem{innercustomthm}{Theorem}
\newenvironment{customthm}[1]
{\renewcommand\theinnercustomthm{#1}\innercustomthm}
{\endinnercustomthm}

\newtheorem{innercustomlemma}{Lemma}
\newenvironment{customlemma}[1]
{\renewcommand\theinnercustomlemma{#1}\innercustomlemma}
{\endinnercustomlemma}

\newtheorem{innercustomprop}{Proposition}
\newenvironment{customprop}[1]
{\renewcommand\theinnercustomprop{#1}\innercustomprop}
{\endinnercustomprop}

\theoremstyle{definition}

\newtheorem{innercustomquestion}{Question}
\newenvironment{customquestion}[1]
{\renewcommand\theinnercustomquestion{#1}\innercustomquestion}
{\endinnercustomquestion}

\allowdisplaybreaks

\renewcommand{\1}{\mathbbm{1}}

\def\L{\mathbf{l}}
\def\R{\mathbf{r}}
\def\B{\mathbf{b}}
\def\1{\mathbf{1}}
\def\p{\mathtt{p}}
\def\q{\mathtt{q}}
\def\r{\mathtt{r}}

\begin{document}

\vspace{-0.5cm}

\title[Clustering in cyclic particle systems on $\mathbb{Z}$]{Clustering in the three and four color \\cyclic particle systems in one dimension}

\author{Eric Foxall}
\address{Eric Foxall, Department of Mathematics, University of Alberta, Edmonton, AB, Canada}
\email{\texttt{e.t.foxall@gmail.com}}

\author{Hanbaek Lyu}
\address{Hanbaek Lyu, Department of Mathematics, The Ohio State University, Columbus, OH, USA}
\email{\texttt{colourgraph@gmail.com}}

\date{\today}


\begin{abstract}
We study the $\kappa$-color cyclic particle system on the one-dimensional integer lattice $\mathbb{Z}$, first introduced by Bramson and Griffeath in \cite{bramson1989flux}. In that paper they show that almost surely, every site changes its color infinitely often if $\kappa\in \{3,4\}$ and only finitely many times if $\kappa\ge 5$. In addition, they conjecture that for $\kappa\in \{3,4\}$ the system clusters, that is, for any pair of sites $x,y$, with probability tending to 1 as $t\to\infty$, $x$ and $y$ have the same color at time $t$. Here we prove that conjecture.
\end{abstract}

\maketitle

\section{Introduction}
\label{Section: Introduction}

Let $G=(V,E)$ be a simple graph and for each vertex (site) $v\in V$, denote by $N(x)$ the set of all neighbours of $x$ in $G$. Fix an integer $\kappa \ge 2$. A \textit{$\kappa$-coloring} on $G$ is a mapping $X:V\rightarrow \mathbb{Z}_{\kappa}=\mathbb{Z}/\kappa \mathbb{Z}$. Let $X_{0}$ be the random initial $\kappa$-coloring on $G$ drawn from the uniform product measure on $(\mathbb{Z}_{\kappa})^{V}$, and as in \cite{bramson1989flux}, for each ordered pair $(x,y)$ such that $\{x,y\} \in E$, define an independent Poisson point process $U(x,y)$ with intensity 1. The \textit{$\kappa$-color cyclic particle system} (CPS) on $G$ is the stochastic process $(X_{t})_{t\ge 0}$ of $\kappa$-colorings on $G$, such that if $t \in U(x,y)$ and $X_{t^-}(y)=X_{t^-}(x)-1$ $(\text{mod $\kappa$} )$ then $X_t(y)=X_{t^-}(x)$. The case $\kappa=3$ can be thought of as a repeated game of ``rock paper scissors'' between agents, in which adjacent pairs of agents play at rate $1$ and the loser adopts the winner's strategy after each interaction. We say the process $X_{t}$ \textit{fluctuates} if every vertex changes its color infinitely often, and \textit{fixates} otherwise. We also say the process \textit{clusters} if 
\begin{equation*}
	\lim_{t\rightarrow \infty} \mathbb{P}(X_{t}(x)= X_{t}(y)) = 1 \quad \text{for all $x,y\in V$}. 
\end{equation*}

The $\kappa$-color cyclic particle system was first introduced by Bramson and Griffeath  in 1989 \cite{bramson1989flux}. The authors studied this process on the one-dimensional integer lattice $V=\mathbb{Z}$ with nearest-neighbour edges, and showed that the process fluctuates if $\kappa\le 4$ and fixates otherwise. In addition, they conjectured that the process clusters if $\kappa\in \{3,4\}$. In this paper, we prove that conjecture.

\begin{customthm}{1}\label{thm:mainthm}
	Let $(X_{t})_{t\ge 0}$ be the $\kappa$-color cyclic particle system trajectory on $\mathbb{Z}$. If $\kappa\in \{3,4\}$, then $X_{t}$ clusters.
\end{customthm} 

\hspace{-0.44cm}We remark that our approach only uses that the distribution of $X_0$ is spatially ergodic and invariant with respect to both spatial reflection and cyclic color shift. Using the ergodic decomposition theorem, it follows that Theorem \ref{thm:mainthm} is true assuming only translation, reflection and color shift invariance. 
We also note that for $\kappa=2$, the CPS reduces to the classical voter model with two opinions, which is known to fluctuate and cluster in one dimension \cite{liggett1997stochastic}. 

\begin{figure*}[h]
	\centering
	\includegraphics[width=0.95 \linewidth]{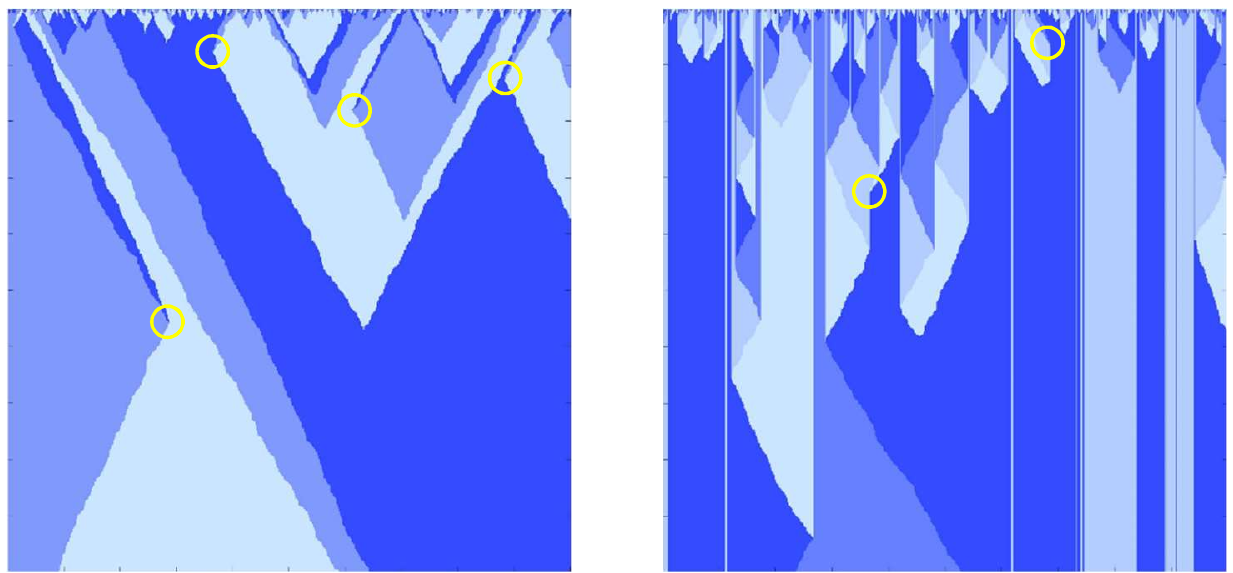}
	\caption{ Simulation of $3$-color (left) and 4-color (right) cyclic particle system on a path of 1000 nodes shown in space-time plots. Time runs downwards and the rows show the colorings at corresponding times. The yellow circles to the left (resp., right) indicate some space-time locations of particle flipping (resp., blockade creation). See Section Section \ref{Section:embedded particle system} for definitions of these terminology. 
	}
	\label{fig:simulation}
\end{figure*}

This paper is organized as follows. In Section \ref{Section:review}, we review some of the related works and state some open questions. In Section \ref{Section:embedded particle system}, we give a precise description of the embedded system of edge particles for the 3- and 4-color CPS on $\mathbb{Z}$, which forms the basis of our discussion. In Section \ref{Section:proof} we prove our main result, Theorem \ref{thm:mainthm}.

\vspace{0.3cm}

\section{Related works and further questions}
\label{Section:review}

In this section, we review some related works and state some open problems.

\subsection{Multitype voter models} The cyclic particle system belongs to the literature of interacting particle systems, many of which can be thought of as opinion exchange dynamics on a social network. The general setup is as follows: let $G=(V,E)$ be a connected simple graph whose nodes and edges represent agents and connection between them, respectively. At each time $t\ge 0$, each agent $v\in V$ may have one of the available opinions $[\kappa]:=\{0,\cdots,\kappa-1\}$. As in the CPS, each ordered pair $(x,y)$ of adjacent sites has a Poisson point process $U(x,y)$ of intensity $1$. Each agent $y$ updates her current opinion according to some transition rule at times $t\in U(x,y)$ for some $x\in N(y)$. This generates a continuous-time Markov chain of opinions $(\eta_{t})_{t\ge 0}$, $\eta_{t}:V\rightarrow [\kappa]$.  

In the classical multitype voter model, each time an agent $x$ updates her opinion at time $t\in U(x,y)$, she simply copies $x$'s opinion, i.e., $\eta_{t}(y)=\eta_{t^{-}}(x)$ (see, e.g., Liggett \cite{liggett1997stochastic}). This rule may be generalized as follows. Let $\Gamma=([\kappa],F)$, $F\subset [\kappa]^{2}$ be a directed graph on the set $[\kappa]$ of available opinions, which we may call the \textit{opinion graph}. The corresponding rule is that at any time $t\in U(x,y)$, $y$ copies $\eta_{t^{-}}(x)$ iff $(\eta_{t}(x),\eta_{t}(y))\in F$. It is natural to view this model as competition dynamics between $\kappa$ distinct species inhabiting the nodes of $G$ with a `general appetite rule' given by $\Gamma$ (see the discussion of Bramson and Griffeath in \cite{bramson1989flux}). 

In the above setting, we recover the classical multitype voter model by taking $\Gamma$ to be the complete graph on $[\kappa]$. The $\kappa$-color cyclic particle system is obtained by letting $\Gamma$ be the directed cycle $\mathbb{Z}_{\kappa}$, where $(i,j)\in F$ iff $i\equiv j+1$ (mod $\kappa$). Recently, Lanchier and Scarlatos \cite{lanchier2014limiting} studied the model on the vertex set $V=\mathbb{Z}$ in the case where $\Gamma$ is symmetric, i.e., $(i,j)\in F$ iff $(j,i)\in F$. The authors give sufficient conditions for fluctuation, clustering and fixation to occur. In particular, they show that the process clusters if $\Gamma$ has radius 1, that is, there exists an opinion $i\in [\kappa]$ which is adjacent to all the other opinions in $\Gamma$ (Theorem 1.1 in \cite{lanchier2014limiting}). In our case, $\Gamma$ is not symmetric and its radius is $\kappa-1$, so the result of Theorem \ref{thm:mainthm} captures a different phenomenon.

\subsection{Edge particle systems and cellular automata} 
\begin{figure*}[h]
	\centering
	\includegraphics[width=1 \linewidth]{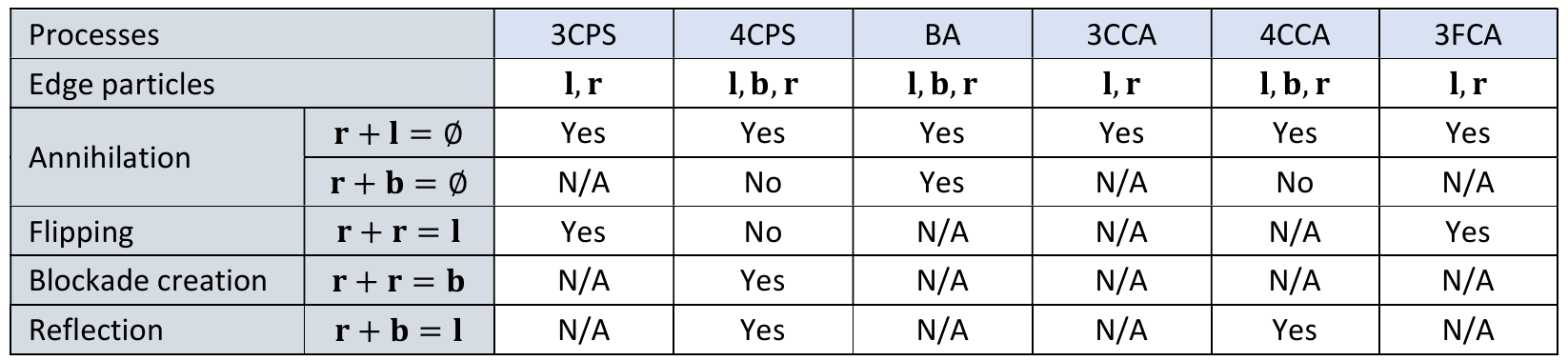}
	\caption{ Collision rules for $\R$ particles in the 3- and 4-color cyclic particle systems, ballistic annihilation, 3- and 4-color cyclic cellular automata, and the 3-color firefly cellular automaton. Rules for $\L$ particles are defined symmetrically. 
	}
	\label{fig:table}
\end{figure*}

A common tool for analyzing interacting particle systems on $\mathbb{Z}$ is to place edge particles at the boundaries between monochromatic regions and follow their time evolution. Here we give an informal description of the dynamics of these particles for the 3- and 4-color CPS, as well as some other models; a formal derivation is provided in Section \ref{Section:embedded particle system}. In the 3-color CPS, between any two monochromatic regions containing at least 2 vertices, the boundary between the two regions can move either left or right (but not both), depending on their respective types. If it can move left (resp. right) we place a left (resp. right) particle, also denoted $\L$ or $\la$ ($\R$ or $\ra$). In the 4-color CPS there is an additional third type, the blockade ($\B$), that is placed at the boundary between regions whose color differs by 2, and does not move so long as it persists. When one particle tries to occupy the same edge as another particle already on that edge, they can be removed and/or change type; see Section \ref{Section:embedded particle system} for a complete description. The collision rules for the 3- and 4-color CPS, and some other similar models, are given in Figure \ref{fig:table}.

Ballistic annihilation (BA) was introduced in the 1990's in the physics community \cite{ben1993decay}. There, particles are released from a Poisson point process (or any translation invariant distribution) on $\mathbb{R}$  with i.i.d. velocities sampled from a bounded distribution $\mu$. Once we place particles and assign their velocities, they move deterministically with a constant velocity, annihilating with one other upon collision. A specific case that has been given a detailed treatment is when $\mu$ is symmetric on $\{-1,0,1\}$ \cite{dygert2016bullet}. In this setup, there are three types of particles ($\L$, $\B$, and $\R$), assigned according to their initial velocities, and they interact by annihilating upon collision (see Figure \ref{fig:table}). See \cite{sidoravicius2017note} for a recent survey on this process.

When $\mu(0)=0$, i.e., when there are no blockades initially, then using a connection with the persistence of an associated random walk or running maximum of comparison surface growth model, it can be shown that particle density decays to zero on the order $t^{-1/2}$ as time $t$ tends to infinity \cite{elskens1985annihilation, krug1988universality, belitsky1995ballistic}. This is in fact the embedded edge particle system structure of the 3-color \textit{cyclic cellular automaton} (CCA), a discrete time version of the cyclic particle system introduced and studied by Fisch \cite{fisch1990cyclic}. Namely, a given $\kappa$-coloring $Y_{t}:V\rightarrow \mathbb{Z}_{\kappa}$ on $G$ updates in discrete time according to the following rule:
\begin{equation*}
	Y_{t+1}(v) = \begin{cases}
		Y_{t}(v)+1 \mod \kappa & \text{if $\exists u\in N(v)$ s.t. $Y_{t}(u)=Y_{t}(v)+1 \mod \kappa$} \\
		Y_{t}(v) & \text{otherwise}.
	\end{cases}
\end{equation*} 
Let $G=\mathbb{Z}$ and draw the initial $\kappa$-coloring $Y_{0}$ uniformly at random. Analogous to the cyclic particle systems, Fisch showed that the $\kappa$-color CCA trajectory $(Y_{t})_{t\ge 0}$ fluctuates for $\kappa\in \{3,4\}$, clusters for $\kappa=3$, and fixates for $\kappa\ge 5$ \cite{fisch1990one, fisch1992clustering}.

\begin{figure*}[h]
	\centering
	\includegraphics[width=1 \linewidth]{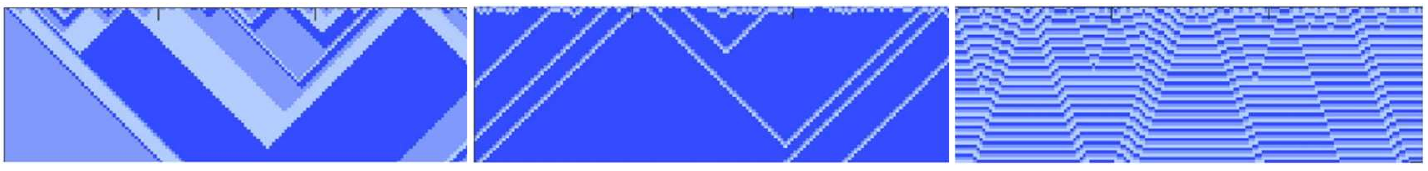}
	\caption{ Simulation of $3$-color CCA (left), GHM (middle), and FCA (right) on a path of 400 nodes for 50 iterations. The top rows are a random $3$-color initial coloring, and a single iteration of the corresponding transition map generates each successive row from top to bottom. Dark blue=0, blue=1, and light blue=2.
	}
	\label{fig:matching}
\end{figure*}

A characteristic collision rule for the 3-color CPS is `flipping', which occurs when an $\L$ or $\R$ particle collides with a particle of the same type (see Figure \ref{fig:simulation} left). While this type of collision does not occur in both BA and CCA due to the synchronous movement of same-typed particles, it does occur in another related cellular automata called the \textit{firefly cellular automata} (FCA). The FCA was introduced by the second author in \cite{lyu2015synchronization} as a discrete model for pulse-coupled oscillators such as fireflies and circadian pacemaker cells. In its embedded particle system, flipping does occur but only for the first finite number of iterations.

For the 4-color CPS and CCA, edges with color difference 2 give rise to blockades. In the 4-color CCA, blockades are never created, and upon collision they reflect $\L$ or $\R$ particles and are removed. A simple ergodic argument could show that this system clusters on $\mathbb{Z}$, when initialized from the uniform product measure. In the 4-color CPS, blockades follow the same reflection rule when hit by a $\L$ or $\R$ particle, but are also created by $\L$-$\L$ or $\R$-$\R$ collisions (see Figure \ref{fig:simulation} right).

\subsection{Further questions: Clustering rate} 

Our proof of Theorem \ref{thm:mainthm} relies on Birkoff's ergodic theorem so it gives no information on the clustering rate, that is, the speed at which particle density decays to zero. This leaves the following question.

\begin{customquestion}{1}\label{question1}
	Let $\kappa\in \{3,4\}$ and let $X_{t}$ be a $\kappa$-color CPS trajectory on $\mathbb{Z}$ started from the uniform product measure. Are there some constants $C,\alpha>0$ depending on $\kappa$, such that 
	\begin{equation*}
	\mathbb{P}(X_{t}(x)\ne X_{t}(x+1)) \sim Ct^{-\alpha}?
	\end{equation*}  
\end{customquestion}

Results of this kind were obtained for the 3-color CCA (Fisch \cite{fisch1992clustering}), $\kappa$-color GHM (Durrett and Steif \cite{durrett1991some} and Fisch and Gravner \cite{fisch1995one}), and recently for the $\kappa$-color FCA (Lyu and Sivakoff \cite{lyu2017persistence, lyu2017synchronization}), for all $\kappa\ge 3$. We shall give a brief sketch of derivation of clustering rates for the 3-color CCA and the 3-color FCA \cite{lyu2017persistence}. This will help the readers to see the difficulty in proving  analogous clustering results for the CPS in Question \ref{question1}.

Fix a 3-color CCA trajectory $(Y_{t})_{t\ge 0}$, started from the uniform product probability measure. We would like to show that 
\begin{equation*}
\mathbb{P}(Y_{t}(x)\ne Y_{t}(x+1)) \sim \sqrt{\frac{2}{3\pi}}t^{-1/2}
\end{equation*}
for all $x\in \mathbb{Z}$. Based on the simulation in Figure \ref{fig:simulation}, it is natural to imagine right- and left-moving particles $\R$ and $\L$ on the boundaries of monochromatic regions given by the initial 3-coloring $Y_{0}:\mathbb{Z}\rightarrow \mathbb{Z}_{3}$, with the particle pointed in the direction of the color that is one value smaller mod 3. These particles move at unit speed in parallel without changing their directions; if two opposing particles ever cross or have to occupy the same edge, then they annihilate each other. 

Thus, if one wants to have an $\R$ particle on the edge $(0,1)$ at time $t$, which occurs with probability the half of $\mathbb{P}(Y_{t}(0)\ne Y_{t}(1))$ by symmetry, then this $\R$ particle must have been on the edge $(-t,-t+1)$ at time $0$ and it must travel distance $t$ without being annihilated by an opposing particle. Thus, at time 0, on the edges of every interval $[-t,s]$ for $-t+1\le s \le t+1$, there must be more $\R$ than $\L$ particles. The converse of this observation is also true. So this running sum of $\R$ versus $\L$ particles starting from the edge $(-t,-t+1)$ to the right must stay above level 1. Now since $Y_{0}$ is drawn from the uniform product measure on $(\mathbb{Z}_{3})^{\mathbb{Z}}$, the running sum of particles at time 0 has i.i.d. centered increments, and the precise asymptotic of it staying above a certain level in the first $t$ steps can be computed by Sparre-Anderson's formula. 

For the 3-color FCA, initial $\R$ and $\L$ particles are assigned similarly. Some of them may flip their directions during the first iteration, and thereafter they follow the similar annihilating dynamics where particles move at speed $1/3$ (see Figure \ref{fig:simulation}). Consequently, one has an $\R$ particle on the edge $(0,1)$ at time $3t+1$ if and only if there is an $\R$ particle on the edge $(-t,-t+1)$ at time $1$ and the running particle sum on the interval $[-t,t+1]$ stays positive. Unlike in the CCA case, increments of the running sum at time 1 have long range correlation. Recently, Lyu and Sivakoff \cite{lyu2017persistence} generalized Sparre-Anderson's formula to handle this type of situation and obtained a precise clustering rate for the 3-color FCA. 

The premise of this kind of running sum argument is that after a certain time $t=t_{0}$, assigned particles do not flip their direction so that particle at a given location at a future time $t>t_{0}$ can be traced back to time $t=t_{0}$. However, flipping and blockade creation may occur arbitrarily many times in the 3- and 4-color CPS so a similar argument does not carry over. We remark that handling blockades is a key difficulty in obtaining the clustering rate for the 4-color CCA (see  \cite{fisch1992clustering}) and rigorously proving fluctuation in the 3-speed BA (see \cite{sidoravicius2017note}).

\vspace{0.3cm}

\section{Embedded particle system}
\label{Section:embedded particle system}

In this section, we give a precise description of the embedded system of edge particles for the 3- and 4-color CPS on $\mathbb{Z}$. Let $E = \{ x+1/2\,:\, x \in \Z \}$ index the set of nearest-neighbour edges on the vertex set $V=\mathbb{Z}$. Given a $\kappa$-coloring $X:\mathbb{Z}\rightarrow \mathbb{Z}_{\kappa}$, $\kappa\in \{3,4\}$, define a map $dX:E\rightarrow \{-1,0,1,\lfloor \kappa/2 \rfloor \}$ by 
\begin{equation*}
dX(e) = X(e+1/2)-X(e-1/2) \,\, (\text{mod $\kappa$}).
\end{equation*} 
Values of $-1,1,2$ and $0$ correspond to $\R,\L,\B$ particles and vacant edges, respectively. $\R$ (right or $\ra$) and $\L$ (left or $\la$) are \emph{directed particles}, and $\B$ is a \textit{blockade}. We call $dX$ the \textit{edge particle configuration} induced by $X$. The process of edge particle configurations $(dX_{t})_{t\ge 0}$ is called the \textit{embedded edge particle system} of $(X_{t})$.\\

\begin{figure*}[h]
	\centering
	\includegraphics[width=0.65 \linewidth]{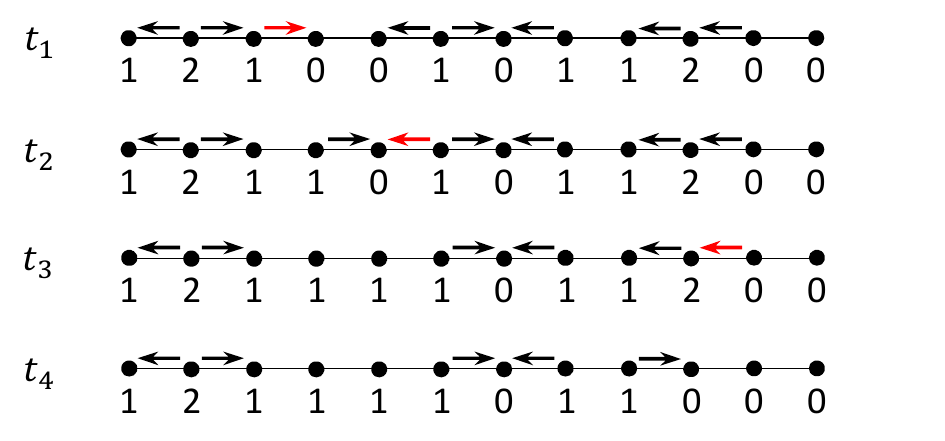}
	\caption{ Example of 3CPS trajectory on a finite path. Numbers indicate colors, Arrows indicate $\L$ and $\R$ particles where red ones are in action in each transition from time $t_{i}$ to $t_{i+1}$. The three transitions show the events `Move', `Annihilate', and `Flip', respectively.     
	}
	\label{fig:3CPS}
\end{figure*}

\begin{figure*}[h]
	\centering
	\includegraphics[width=0.65 \linewidth]{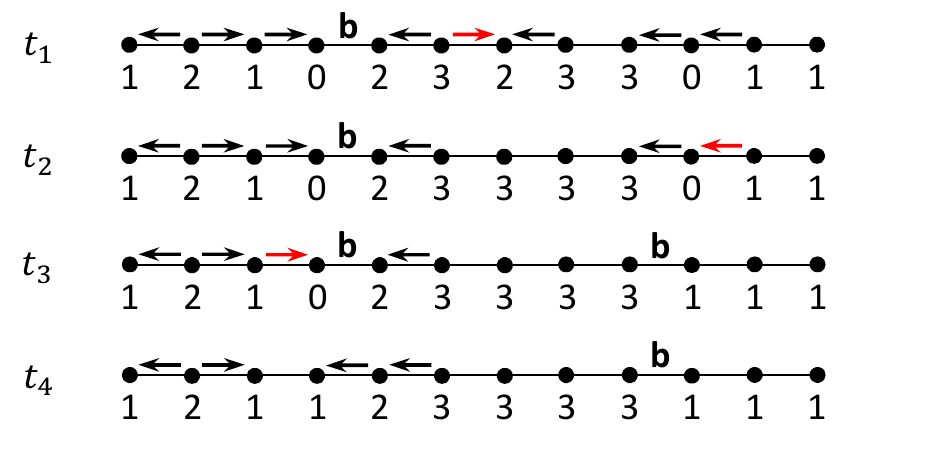}
	\caption{ Example of 4CPS trajectory on a finite path. Numbers indicate colors, $\B$'s indicate blockades, arrows indicate $\L$ and $\R$ particles where red ones are in action in each transition from time $t_{i}$ to $t_{i+1}$. The three transitions show the events  `Annihilate', `Blockade creation', and `Reflect', respectively.    
	}
	\label{fig:4CPS}
\end{figure*}

To describe the update rules for $(dX_t)_{t \ge 0}$ we use the intuitive language of particle motion, collision, and change of type. For example, if we say that an $\R$ particle moves from $e$ to the vacant edge $e+1$ at time $t$, it means that $dX_{t^-}(e)=-1$, $dX_{t^-}(e+1)=0$, $dX_t(e)=0$ and $dX_t(e+1)=-1$. Figures \ref{fig:3CPS} and \ref{fig:4CPS} help to make the updates clear for 3 and 4 colors, respectively. Recall from the introduction that to each edge we associate two Poisson point processes with intensity 1, one pertaining to each direction along that edge. We will denote these by $U^+(e)$ (from $e-1/2$ to $e+1/2$) and $U^-(e)$ (vice-versa).

\subsection{Evolution of edge particles for $\kappa=3$}\label{subsection:particle_3color} Suppose that $dX_{t^-}(e)=-1$ (there is an $\R$ particle on $e$ just before time $t$) and $t \in U^+(e)$. Then, it is easy to check that the following occurs:

\begin{description}
	\item{(i)} (Move) If edge $e+1$ is vacant at time $t^-$ then the $\R$ particle on $e$ moves to $e+1$.
	\item{(ii)} (Annihilate) If the edge $e+1$ has an $\L$ particle at time $t^-$, then both the $\R$ on $e$ and $\L$ on $e+1$ are removed from the system. 
	\item{(iii)} (Flip) If the edge $e+1$ has an $\R$ particle at time $t^-$, then this $\R$ is removed from the system and the $\R$ on $e$ at time $t$ becomes an $\L$ on $e+1$. 
\end{description}

Behavior of $\L$ particles is the same except in mirror image and using $U^-(e)$ instead of $U^+(e)$. The second and third options are collectively referred to as a \emph{collision}. It is important to note that two particles are involved in every collision, and at least one is removed.

\subsection{Evolution of edge particles for $\kappa=4$} \label{subsection:particle_4color} A blockade remains immobile for as long as it persists, since sites whose colors differ by 2 do not interact. Suppose that edge $e$ has an $\R$ particle at time $t^-$ and $t \in U^+(e)$. Then, it is easy to check that the following occurs:

\begin{description}
	\item{(i)} (Move) If edge $e+1$ is vacant at time $t^-$, the $\R$ particle on $e$ moves to $e+1$.
	\item{(ii)} (Annihilate) If edge $e+1$ has an $\L$ particle at time $t^-$, then both the $\R$ on $e$ and $\L$ on $e+1$ are removed from the system. 
	\item{(iii)} (Blockade creation) If the edge $e+1$ has an $\R$ particle at time $t^-$, this $\R$ is removed from the system and the $\R$ on $e+1$ at time $t$ becomes a blockade $\B$ on $e+1$. 
	\item{(iv)} (Reflect) If edge $e+1$ has a blockade $\B$ at time $t^-$, this $\B$ is removed from the system and the $\R$ on $e$ at time $t$ becomes an $\L$ on $e+1$. 
\end{description}     

Behavior of $\L$ particles is the same except in mirror image. The last three options are collectively referred to as a \emph{collision}. As in the 3-color case, every collision involves two particles, and at least one particle is removed.

\vspace{0.3cm}
\section{Proof of Theorem \ref{thm:mainthm}} 
\label{Section:proof}

Let $(X_{t})_{t\ge 0}$ be the $\kappa$-color cyclic particle system trajectory started from the uniform product measure, for $\kappa\in \{3,4\}$. Since the distribution of $X_0$, and the update rules, are translation and color shift invariant, for any $t\ge 0$ and $a \in \{1,2\}$ the probabilities $\mathbb{P}(|dX_{t}(e)|=a)$ for $a\in \{1,2\}$ do not depend on $e\in E$. For each $t\ge 0$, define the respective density of directed particles and blockades by
\begin{equation*}
\mathtt{p}(t) = \mathbb{P}(|dX_{t}(e)|=1) \quad\text{and}\quad \mathtt{q}(t) = \mathbb{P}(|dX_{t}(e)|=2),
\end{equation*}
and let $\r(t) = \p(t)+\q(t)$ denote the total density of particles. To show clustering occurs it is enough to show that $\lim_{t\to\infty}\r(t)=0$, since by a union bound, for $x<y$
$$P(X_t(x) \ne X_t(y)) \le P(dX_t(e) \ne 0 \ \text{for some} \ e \in [x,y]) \le (y-x)\r(t).$$
The proof of Theorem \ref{thm:mainthm} follows from two lemmas. Recall that, as shown in \cite{bramson1989flux}, for $\kappa\in \{3,4\}$ the $\kappa$-color CPS trajectory on $\mathbb{Z}$ fluctuates. Using fluctuation, together with the mass-transport principle, we will show the density of directed particles is at least as large as the density of blockades (which is trivial for $\kappa=3$ as $\mathtt{q}(t)\equiv 0$). This is stated in the following lemma. 

\begin{customlemma}{4.1}\label{lemma:fixation}
	Fix $\kappa\in \{3,4\}$. Then $\mathtt{p}(t)\ge \mathtt{q}(t)$ for all $t\ge 0$. 
\end{customlemma}

Thus it suffices to show that $\p(t) \to 0$ as $t\rightarrow \infty$ in order to establish clustering. The following lemma, Lemma \ref{lemma:key}, shows that the total density of particles does not increase, and that a constant fraction of directed particles is subtracted from the total in any long enough time window:

\begin{customlemma}{4.2}\label{lemma:key}
	Fix $\kappa\in \{3,4\}$. Then, $r(s) \le r(t)$ for $s\ge t \ge 0$, and for each $t\ge 0$, there is $s=s(t)\ge t$ such that
	\begin{equation*}
		\r(s) \le \r(t) - \p(t)/4.
	\end{equation*}
\end{customlemma}

We prove Theorem \ref{thm:mainthm} based on the two lemmas above. 

\vspace{0.1cm}
\hspace{-0.42cm}\textbf{Proof of Theorem \ref{thm:mainthm}.} By Lemma \ref{lemma:fixation}, $\mathtt{r}(t) \le 2\mathtt{p}(t)$ for every $t\ge 0$ and by Lemma \ref{lemma:key}, $\r(t)$ is non-increasing in $t$, so it suffices to find a sequence $t_0\le t_1 \le \dots $ such that $p(t_n) \to 0$ as $n\to\infty$. Let $t_0=0$ and iterating Lemma \ref{lemma:key}, let $t_{n+1}=s(t_{n})$ for all $n\ge 0$. Then
\begin{equation*}
\mathtt{r}(t_n) \le \mathtt{r}(0) - \frac{1}{4}\left( \mathtt{p}(0)+\mathtt{p}(t_1) + \cdots+\mathtt{p}(t_{n-1}) \right)
\end{equation*}
Since $\mathtt{r}(t_n) \ge 0$ for all $n$, letting $n\to\infty$ we see that $\sum_n\mathtt{p}(t_{n})$ converges. In particular, $\mathtt{p}(t_{n})\rightarrow 0$ as $n\rightarrow \infty$. $\blacksquare$ 

\vspace{0.2cm}
The rest of this section is devoted to the proof of Lemmas \ref{lemma:fixation} and \ref{lemma:key}. As mentioned earlier, the former follows from fluctuation and the mass transport principle.  

\vspace{0.1cm}
\hspace{-0.42cm}\textbf{Proof of Lemma \ref{lemma:fixation}.} Since $\mathtt{q}(t)\equiv 0$ for $\kappa=3$, we may assume $\kappa=4$. For each $t'>t\ge 0$, define the density of blockades at time $t$ surviving through time $t'>t$:
	\begin{equation*}
	\mathtt{q}(t,t') = \mathbb{P}( \ dX_s(e)=2 \ \forall s \in [t,t'] \ ).
	\end{equation*} 
	Note the above probability does not depend on $e\in E$. Moreover, $\mathtt{q}(t,t')$ is non-increasing in $t'$. For each $t\ge 0$, denote $\delta_{t}:=\lim_{t'\rightarrow \infty} \mathtt{q}(t,t')$.
	
	First we claim that $\delta_{t}=0$ for all $t\ge 0$. To see this, suppose $\delta_{t}>0$ for some $t\ge 0$. By ergodicity with respect to translation on $\Z$, an asymptotic fraction $\delta_t$ of edges have blockades that are unchanged on the time interval $[t,\infty)$, so in particular there is at least one such edge; denote it $e$. Since the state at adjacent sites does not change simultaneously, it follows that $(X_s(e-1/2),X_s(e+1/2))$ is constant for $s \in [t,\infty)$. But this contradicts fluctuation.
	
	Next, for $t\ge 0$ and $x,y\in \mathbb{Z}$, define the indicator variable 
	\begin{equation*}
	Z_{t}(x,y) = \mathbf{1}\left\{   
	\begin{matrix}
	\text{ $|dX_{t}(x+1/2)|=1$, $dX_{t}(y+1/2)=2$ and } \\
	\text{ the directed particle on $x+1/2$ at time $t$} \\
	\text{ collides with  the $\B$ on $y+1/2$ }
	\end{matrix}
	\right\} .
	\end{equation*}  
	We apply an elementary form of the `mass-transport principle' for these indicator variables. Namely, by linearity of expectation and translation invariance of the process $(X_{t})_{t\ge 0}$ on $\Z$, we obtain
	\begin{equation*}
	\mathbb{E}\left[ \sum_{y\in \mathbb{Z}} Z_t(0,y) \right] =  \sum_{y\in \mathbb{Z}} \mathbb{E} Z_t(0,y)  = \sum_{y\in \mathbb{Z}} \mathbb{E} Z_t(-y,0) =  \sum_{y\in \mathbb{Z}} \mathbb{E} Z_t(y,0) = \mathbb{E}\left[ \sum_{y\in \mathbb{Z}} Z_t(y,0) \right].
	\end{equation*}
	Note that the third equality above uses the fact that sum over all $y$'s equals to sum over all $-y$'s. On the one hand, we have  
	\begin{eqnarray*}
	\mathbb{E}\left[ \sum_{y\in \mathbb{Z}} Z_t(0,y) \right] 
	& = & \mathbb{P}\left(  
	\begin{matrix}
	\text{$\exists$ a $\R$ or $\L$ on the edge $1/2$ at time $t$} \\
	\text{that eventually collides with a $\B$} 
	\end{matrix}
	\right) \\
	&\le& \mathbb{P}(|dX_{t}(1/2)|=1)=\mathtt{p}(t).
	\end{eqnarray*}
	On the other hand, by the claim, we have
	\begin{eqnarray*}
		\mathbb{E}\left[ \sum_{y\in \mathbb{Z}} Z_t(y,0) \right] &=& \mathbb{P}\left(  
		\begin{matrix}
		\text{$\exists$ a $\B$ on the edge $1/2$ at time $t$} \\
		\text{that is eventually removed} 
		\end{matrix}
		\right) \\
		&=&\mathbb{P}(dX_{t}(1/2)=2)- \mathbb{P}(\text{$dX_{t'}(1/2)=2$ for all $t'\ge t$})  \\
		&=& \mathbb{P}(dX_{t}(1/2)=2) = \mathtt{q}(t). 
	\end{eqnarray*}
	This shows the assertion. $\blacksquare$

\vspace{0.2cm}

Now it remains to prove Lemma \ref{lemma:key}. Our general strategy is to show that at any given time $t$, a positive density of remaining directed particles are involved in some type of collision in some time window $[t,s]$. For this we introduce some terminology. 
Let $\xi:E\to \{-1,0,1,2\}$ be an edge configuration. For each interval $[a,b]$ with $a,b \in \mathbb{Z}$, a set $\mathcal{M}$ of ordered pairs $(e,e')$ of edges  in $[a,b]$ is called a \textit{$\xi$-matching} if the following two conditions are satisfied: 
\begin{description}
	\item{(i)} 	$(e,e') \in \mathcal{M}  $ implies  \\
	 $e<e'$ and $\xi(e)=-1$ and $\xi(e')=1$;
	\item{(ii)} No two elements of $\mathcal{M}$ share an edge.
\end{description}
Point (i) means that the $\R$ and $\L$ particles on the edges $e$ and $e'$, are matched by $\mathcal{M}$. Point (ii) just means no edges are counted twice in a matching. Since $\mathcal{M}$ consists of pairs of edges, the number of particles matched by $\mathcal{M}$ is $2|\mathcal{M}|$.

It is convenient to introduce virtual particles in order to control the movement of edge particles. Virtual particles move in the same way as real particles, except that they do not interact with other particles, or change type. A virtual $\R$ particle evolves in time as follows:
\begin{equation*}
v_{t}^{\R} - v_{t-}^{\R} = \mathbf{1}( t \in U^+(v_{t^-}^\R) ).
\end{equation*}
A virtual $\L$ particle is defined similarly, except in mirror image and using $U^-$ instead of $U^+$. For $t\ge 0$ and $a',b' \in E$ with $a'<b'$, define the collision time
\begin{equation}\label{eq:def_hittingtime_virtual}
\tau_{t}(a',b') = \inf\{ s\ge t \,:\, v^{\R}_{s} = v^{\L}_{s}   \},
\end{equation}   
where $v^{\R}_{s}$ and $v^{\L}_{s}$ are trajectories of virtual $\R$ and $\L$ particles with $v^{\R}_{t}=a'$ and $v^{\L}_{t}=b'$.

For each $e\in E$ and $s\ge t\ge 0$, define the following indicator variable
\begin{equation*}
\mathtt{col}_{t}^{s}(e) = \mathbf{1}\left\{  
\begin{matrix}
\text{$\exists$ a directed particle on the edge $e$ at time $t$} \\
\text{that collides with some other particle by time $s$ } 
\end{matrix}
\right\}.
\end{equation*} 
The following proposition shows that a large matching guarantees a large number of collisions. 

\begin{customprop}{4.3}
	\label{prop:anni}
	Fix $\kappa\in \{3,4\}$ and an interval $[a,b]$ with $a,b \in \Z$, and let $a'=a+1/2$ and $b'=b-1/2$. Let $\xi$ be a $\kappa$-color edge configuration on $\mathbb{Z}$ and $\mathcal{M}$ be a $\xi$-matching on $[a,b]$. Let $(X_{t})_{t\ge 0}$ be the $\kappa$-color cyclic particle system trajectories with $dX_{0}=\xi$. Then 
	\begin{equation*}
		\sum_{a<e<b} \mathtt{col}_{0}^{\tau_{0}(a',b')}(e) \ge |\mathcal{M}|.
	\end{equation*}
\end{customprop}

\begin{proof}
	It suffices to show that for each $(e,e')\in \mathcal{M}$, one of the two particles initially on $e$ and $e'$ collides with some other particle. Let $u^\R_t,u^\L_t$ denote the locations of these particles, with $u^\R_0=e$ and $u^\L_0=e'$, assigning a coffin state to the location of either, if it is removed. Place virtual $\R$ and $\L$ particles on sites $\ell$ and $r$ at time $0$, denoting their locations at time $t$ by $v_{t}^{\R}$ and $v^{\L}_{t}$, respectively. Suppose, for contradiction, the real $\R$ and $\L$ particles have not collided with any particles by time $\tau_0(a',b')$. Since neither has collided, neither has changed type, and since, in addition, real and virtual $\R$ particles move in response to the same point processes $U^+$ (similarly for $\L$ particles, with $U^-$) and $v^\R_0 \le u^\R_0 < u^\L_0 \le v^\R_0$, almost surely
$$v^\R_t \le u^\R_t < u^\L_t \le v^\L_t \quad \text{for all} \quad t \le \tau_0(a',b').$$
Since $v^\R_{\tau_0(a',b')} = v^\L_{\tau_0(a',b')}$, this contradiction shows the (real) particles must collide with each other by time $\tau_{0}(a',b')$.
\end{proof}

Next we show how to obtain a large matching, given some information about particle counts. Fix $\xi:E\rightarrow \{-1,0,1,2\}$ and let
\begin{equation*}
R(\xi, x) = \sum_{y=0}^{x-1} \1(\xi(y+1/2)=-1) \quad \text{and} \quad L(\xi, x) = \sum_{y=0}^{x-1} \1(\xi(y+1/2)=1)
\end{equation*}
be respectively the number of (real) right, left particles in $\xi$ in the interval $[0,x]$. Define the running sum and particle count by
\begin{equation*}
S(\xi,x) = R(\xi,x) - L(\xi,x) \quad \text{and} \quad C(\xi,x) = R(\xi,x) + L(\xi,x).
\end{equation*} 
Also define the running minimum defined by 
\begin{equation*}
m(\xi,x) = \min_{0\le y < x} S(\xi,y).
\end{equation*}

\begin{customprop}{4.4}\label{prop:matching}
	Let $\xi$, $S$, $C$, $\mathcal{M}$ be as above. Then for any integer $N >0$, there exists a $\xi$-matching $\mathcal{M}$ on the interval $[0,N]\subset \mathbb{Z}$ such that
	\begin{equation}\label{eq:matching_size}
	2|\mathcal{M}| = C(\xi,N) - ( \,2 \,|m(\xi,N)| + S(\xi,N) \, ).
	\end{equation}
\end{customprop}

\begin{proof}
	We first define a matching $\mathcal{M}$ and show that it has the correct size as asserted. Let $s:[0,\infty)\rightarrow \mathbb{R}$ be the linear interpolation of the lattice path $S(\xi,\, \cdot \,):\mathbb{N}_{0}\rightarrow \mathbb{Z}$. For $0\le x<y\le N$ let $(x+1/2,y+1/2) \in \mathcal{M}$ if $s(x)=s(y+1)$ and $s(t)>s(x)$ for $x<t<y+1$. Observe that $(x+1/2,y+1/2)\in \mathcal{M}$ implies $\xi(x+1/2)=-1$ and $\xi(y+1/2)=1$. Pictorially, this matching $\mathcal{M}$ can be obtained as follows. Let $\Gamma$ be the graph of $s$ over the interval $[0,N]$. Flip it upside-down to obtain $-\Gamma$, then fill it with water, letting the water spill over the sides at the endpoints $0$ and $N$. The water will pool in various basins; $\mathcal{M}$ matches pairs of points that lie in the same basin at the same half-integer depth.

	\begin{figure*}[h]
		\centering
		\includegraphics[width=0.95 \linewidth]{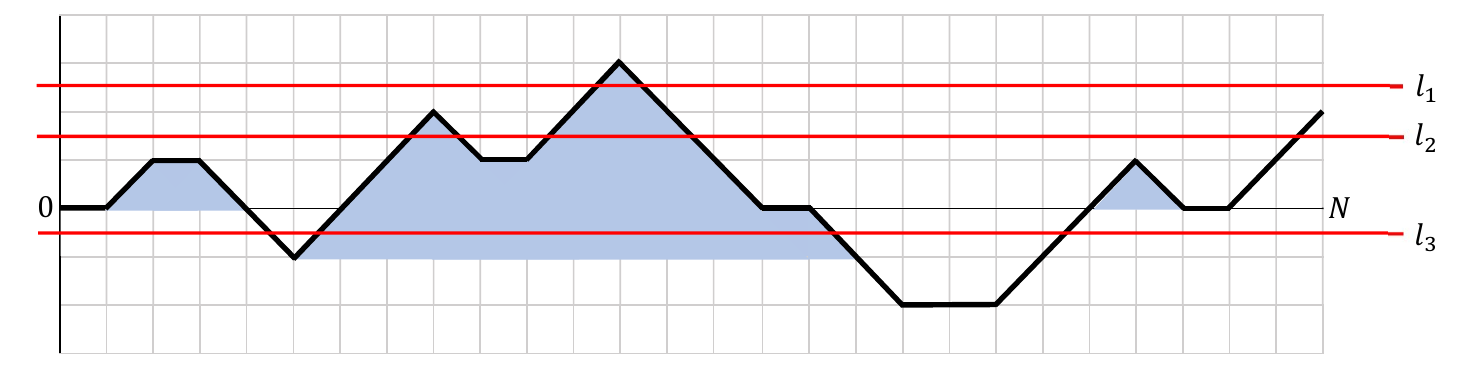}
		\caption{ Linear interpolation $s:[0,N]\rightarrow \mathbb{R}$ of the running sum $S(\xi,x)$ is shown in black. Blue shaded regions depict upside-down image of remaining water after filling the graph of $-s$ with water to the top. Red horizontal lines at half-integer heights show matched and unmatched $\L$ or $\R$ particles.  
		}
		\label{fig:matching}
	\end{figure*}

	Next we show \eqref{eq:matching_size} for the above matching $\mathcal{M}$. To do so, we partition $\R$ and $\L$ particles in $[0,N]$ according to their height on $s$ and match them. For brevity, we introduce the following notation $m^{*}=\max_{[0,N]}s(t)$ and $m_{*}=\min_{[0,N]}s(t)$. For each half-integer $h\in \mathbb{Z}+1/2$ with $m_{*}<h<m^{*}$, the horizontal line $y=h$ intersects with the graph of $s$ over $[0,N]$ in at least one point, and the slope of $s$ alternates between $\pm 1$ at the intersection points, going from left to right. Call such an intersection `upcrossing' (resp., `downcrossing') at height $h$ if $s$ has slope $1$ (resp., $-1$) at that point. Then by construction, $\mathcal{M}$ matches each upcrossing at height $h$ to the next (leftmost to the right) downcrossing at the same height. Note that an intersection at height $h$ is left unmatched by $\mathcal{M}$ iff it is either the leftmost intersection and is a downcrossing, or is the rightmost intersection and is an upcrossing. (See Figure \ref{fig:matching})
	
	Now we consider three cases.  First, let $0\lor s(N)<h<m^{*}$. Then the first and last intersections at height $h$ are up- and down-crossings, respectively, so every intersection is matched by $\mathcal{M}$ (line $l_{1}$ in Figure \ref{fig:matching}). Second, if $0\land s(N)<h< 0\lor s(N)$, then the first and last intersections at height $h$ are both up- or down-crossings, so exactly one intersection at each height $h$ is unmatched (line $l_{2}$ in Figure \ref{fig:matching}). Lastly, if $m_{*}<h< 0\land s(N)$, then the first and last intersections at height $h$ are down- and up-crossings, respectively, so exactly two intersections at each height $h$ are unmatched (line $l_{3}$ in Figure \ref{fig:matching}). Considering what happens in each case, the total number of unmatched particles is equal to $s(N)+2|m_{*}|$ when $s(N)\ge 0$ and equal to $|s(N)|+2(|m_{*}|-|s(N)|) = -s(N) + 2(|m_{*}|+s(N)) = s(N)+ 2|m_{*}|$ when $s(N)<0$. This shows the equality \eqref{eq:matching_size}.        
\end{proof}

Following is an easy observation we will need in further discussions.

\begin{customprop}{4.5}\label{prop:sup}
	Let $(x_{k})_{k\ge 0}$ be a sequence of real numbers such that $x_{k}/k\rightarrow 0$ as $k\rightarrow \infty$. Then $k^{-1}\sup_{1\le i \le k} x_{i}\rightarrow 0$ as $k\rightarrow \infty$. 
\end{customprop}

\begin{proof}
	Omitted.
\end{proof}

Now we show that a positive fraction of directed particles are involved in some collision in any long enough time window. 

\begin{customprop}{4.6}\label{prop:collision_density}
		Fix $t,\epsilon>0$. Then there exists a finite $T=T(t,\eps)\ge 0$ such that for all $e \in E$, 
		\begin{equation*}
		\mathbb{P}(\mathtt{col}_{t}^{t+T}(e)=1) \ge (1-\epsilon)^{3}\mathtt{p}(t)/2.
		\end{equation*}
\end{customprop}

\begin{proof}
	 Note that the distribution of $(dX_{0}(x))_{x\in \mathbb{Z}}$, as well as the Poisson processes on the edges, is reflection invariant and spatially ergodic. Since $(dX_{t})_{t\ge 0}$ is a deterministic function of these variables, the sequence $\{(dX_t(x))_{t\ge 0} \colon x \in \Z\}$ has the same properties. Hence by Birkoff's ergodic theorem, for any fixed $L\in \mathbb{N}$,
	\begin{equation}\label{eq:density_collision_Birkoff}
		\mathbb{P}(\mathtt{col}_{t}^{s}(e)=1) = \frac{1}{L}\mathbb{E} \left[ \sum_{0<e'<L}  \mathtt{col}_{t}^{s}(e') \right] = \lim_{k\rightarrow \infty}\frac{1}{kL} \left[  \sum_{0<e'<kL}  \mathtt{col}_{t}^{s}(e') \right].
	\end{equation}
	Note that $\sum_{0<e'<L}  \mathtt{col}_{t}^{s}(e')$ equals the number of directed particles in the interval $[0,L]$ at time $t$ that collide with some particle by time $s$.

For each $j\ge 0$, let $I_{j}$ be the interval $[jL,(j+1)L]$. Let
$$\tau_{t}^{j}=\tau_{t}\big(jL+1/2, (j+1)L-1/2\big)-t$$
be the amount of time until collision of the virtual $\R$ and $\L$ particles placed at edges $jL+1/2$ and $(j+1)L-1/2$ at time $t$, respectively (see \eqref{eq:def_hittingtime_virtual}). Note that the $\{\tau_{t}^{j}\}_{j \in \Z}$ are i.i.d., and that since the $\R$ or $\L$ particles each move independently at rate $1$, the distribution of $\tau_t^j$ is the sum of $L$ i.i.d. exponential(2) random variables. By Markov's inequality, 
	\begin{equation*}
		\mathbb{P}(\tau_{t}^{j} \ge  L/2\eps) \le  \eps
	\end{equation*}
	for each $j\ge 0$. For each interval $I$ with integer endpoints, denote by $\mathcal{M}_{t}(I)$ the $dX_{t}$-matching of maximum size on the interval $I$. Define the following events: 
	\begin{equation*}
		A_{j} = \left\{ \tau_{t}^{j} < L/2\eps \right\}, \quad B_{j} = \left\{ 2|\mathcal{M}_{t}(I_{j})|>(1-\eps) \mathtt{p}(t)L \right\}.
	\end{equation*}   
	Since $\mathtt{col}_{t}^{t'}(e)$ is non-decreasing in $t'$ for each $e$, Proposition \ref{prop:anni} yields 
	\begin{equation}\label{eq:num_collisions_containment}
		\left\{ \sum_{e\in I_{j}} \mathtt{col}_{t}^{t+L/2\eps}(e) > (1-\eps)\mathtt{p}(t) L/2 \right\} \subseteq A_{j}\cap B_{j}.
	\end{equation}
	Since the event times $U^{\pm}(e) \cap (t,\infty)$, $e \in E$, are independent from the time $t$ configuration $X_{t}$, for each $j\ge 0$, $A_{j}$ and $B_{j}$ are independent. Using this and ergodicity,
	\begin{equation}\label{eq:indicator_frequency}
		\lim_{k\rightarrow \infty} \frac{1}{k} \sum_{j=0}^{k-1} \mathbf{1}(A_{j}\cap B_{j}) = \mathbb{P}(A_{0}) \mathbb{P}(B_{0}) \ge (1-\eps) \mathbb{P}(B_{0}).
	\end{equation}  
	Combining \eqref{eq:num_collisions_containment} and \eqref{eq:indicator_frequency}, for any fixed $L\in \mathbb{N}$,
	\begin{equation}\label{eq:indicator_frequency2}
	\lim_{k\rightarrow \infty} \frac{1}{kL} \sum_{0<e<kL} \mathtt{col}_{t}^{t+L/2\eps} (e)  \ge (1-\eps)^{2} \mathtt{p}(t) \mathbb{P}(B_{0})/2.
	\end{equation}  
	
	To finish, note that since the time $t$ particle configuration $\{(dX_t(x))_{t\ge 0} \colon x \in \Z\}$ is translation and reflection invariant as well as ergodic,
	\begin{equation*}
	\lim_{x\to\infty}x^{-1}R(dX_t,x) = \lim_{x\to\infty}x^{-1}L(dX_t,x) = \mathtt{p}(t)/2
	\end{equation*}
	and hence 
	\begin{equation*}
	\lim_{x\to\infty}x^{-1}S(dX_t,x) = 0, \quad \lim_{x\to\infty}x^{-1}C(dX_t,x) = \mathtt{p}(t) .
	\end{equation*}
	Applying Proposition \ref{prop:sup} to the sequence $x_{k}=-m(dX_{t},k)$,
	\begin{equation*}
	\lim_{x\rightarrow \infty} x^{-1} m(dX_{t},x) =0.
	\end{equation*} 
	Thus by Proposition \ref{prop:matching}, there exists $L_{1}\in \mathbb{N}$ such that $\mathbb{P}(B_{0})>1-\eps$ for all $L\ge L_{1}=L_{1}(t,\eps)$. Let $T(t,\eps)=L_1/2\eps$. Then the assertion follows from \eqref{eq:density_collision_Birkoff} and \eqref{eq:indicator_frequency2}.	
\end{proof}

Lastly, we finish our discussion by proving Lemma \ref{lemma:key}.

\vspace{0.1cm}
\hspace{-0.42cm}\textbf{Proof of Lemma \ref{lemma:key}.}  
It suffices to show that for $s\ge t \ge 0$,
\begin{equation*}\label{eq:q}
\mathtt{r}(s) \le \mathtt{r}(t) - \mathbb{P}(\mathtt{col}_{t}^{s}(e)=1).
\end{equation*}
To see why, note the above implies that $r(t)$ is non-increasing in $t$, which is the first assertion of the lemma. To obtain the second assertion, apply Proposition \ref{prop:collision_density} with $s=t+T(t,\eps)$ and $\epsilon>0$ small enough that $(1-\eps)^3\ge 1/2$.

Given an interval $J=[-N,N]\subset \mathbb{Z}$ and $s\ge 0$, let 
$$r(J,s) = \sum_{-N<e<N}\1(|dX_t(e)| \ne 0)$$
denote the number of edge particles in $J$ at time $s$. Let $t'=s-t$ and define the virtual trajectory $v^{\L}_{t'}$ for $0 \le t' \le s$ by letting $v^{\L}_0=-N+1/2$ and following the rule $v^{\L}_{t'}=e$ if $v^{\L}_{t'^-}=e+1$ and $s-t' \in U^+(e)$. Similarly, define $v^{\R}_{t'}$ for $0\le t' \le s$ by $v^{\R}_0=N-1/2$ and $v^{\R}_{t'}=e$ if $v^{\R}_{t'^-}=e-1$ and $s-t' \in U^-(e)$. For $0 \le t \le s$, define $J_{t} = [v^{\R}_{s-t},v^{\L}_{s-t}]$. Since virtual trajectories move away from each other each at rate 1,
\begin{equation}\label{eq:length_J_t}
	|J_{t}| \overset{d}{=} 2N + \text{Poi}(2(s-t)).
\end{equation}   
Moreover, since particles cannot be created (they are only removed or change type), and cannot cross over the virtual trajectories from the outside to the inside of $J_t$, every particle in $J=J_s$ at time $s$ must be in $J_t$ at every time $t<s$. Since every collision kills at least one particle, it follows that
\begin{eqnarray*}
	\mathtt{r}(J,s) &\le& \# \,\, \text{of particles in $J_{t}$ at time $t$ that survive up to time $s$} \\
	&\le & \mathtt{r}(J_{t},t) - \sum_{v^{\R}_{t}\le e \le v^{\L}_{t}} \mathtt{col}_{t}^{s}(e)\\
	&\le & \mathtt{r}(J_{t},t) - \sum_{-N < e < N} \mathtt{col}_{t}^{s}(e)
\end{eqnarray*}
for every $t<s$. Since $s,t$ are fixed, as $N\to\infty$, $|J_t|/2N \to 1$ almost surely from \eqref{eq:length_J_t}. Using this and spatial ergodicity of $(X_t)_{t \ge 0}$ applied to the functions $\r$ and $\mathtt{col}$,
\begin{eqnarray*}
	\mathtt{r}(s) = \lim_{N\rightarrow \infty} \frac{\mathtt{r}(J,s)}{2N}   &\le&  \lim_{N\rightarrow \infty} \frac{\mathtt{r}(J_{t},t)}{|J_{t}|} \frac{|J_{t}|}{2N} - \lim_{N\rightarrow\infty} \frac{1}{2N} \sum_{-N< e < N} \mathtt{col}_{t}^{s}(e) \\
	&=& \mathtt{r}(t) - \mathbb{P}(\mathtt{col}_{t}^{s}(e)=1),
\end{eqnarray*}
as desired. $\blacksquare$

\vspace{0.5cm}

\small{
	\bibliographystyle{plain}
	\bibliography{mybib}
}

\end{document}